\batchmode
\documentclass[11pt]{amsart}
 \usepackage{amssymb}
\usepackage{amscd} 
\begin{document}
    \newtheorem{Theorem}{Theorem}[section]
    \newtheorem{Proposition}[Theorem]{Proposition}
    \newtheorem{Lemma}[Theorem]{Lemma}
    \newtheorem{Corollary}[Theorem]{Corollary}
   \newcommand{\ra}{\rightarrow}   
   \newcommand{\ul}{\underline}
   \newcommand{\varp}{\varphi} 
   \newcommand{\al}{\alpha}
   \newcommand{\bt}{\beta}
   \newcommand{\gmm}{\gamma}
   \newcommand{\ve}{\varepsilon}
    \newcommand{\aut}{{\rm Aut}} 
     \newcommand{\Z}{\boldsymbol{Z}} 
     \newcommand{\C}{\boldsymbol{C}} 
     \newcommand{\R}{\boldsymbol{R}} 
     \newcommand{\B}{\boldsymbol{P}} 
     \newcommand{\Q}{\boldsymbol{Q}} 
\def\labelenumi{\theenumi)}

\title{Automorphisms of parabolic Inoue surfaces} 

\author{A. Fujiki} 

\maketitle 

\vspace{5 mm} 
\begin{abstract} 
We determine explicitly the structure of the automorphism group of a parabolic Inoue surface. 
We also describe the quotients of the surface by typical cyclic subgroups of 
the automorphism group. 
\end{abstract}

\section{Statement of Results} 
In this note we determine the automorphism group Aut\,$S$ 
of a parabolic Inoue surface $S$.  
The corresponding result for a hyperbolic Inoue surface 
(Inoue-Hirzebruch surface) 
was obtained by Pinkham \cite{pik} more than twenty years ago.  
But the far easier case of parabolic Inoue surfaces does not seem explicit 
in the literature.  

Any parabolic Inoue surface S with second betti number $m >0$ 
contains a unique smooth elliptic curve E and 
a cyclic of rational curves  $C=C_1+\cdots +C_m$, 
where the self-intersection number  $E^2=-m$, and  $C_i$ are nonsingular rational curves 
with $C_i^2=-2$ when $m\geq 2$, while when 
$m=1, C=C_1$ is a rational curve with a single node with $C^2=0$. 
$S$ contains no other curves and hence 
any automorphism of $S$ leaves $C$ and $E$ invariant. 

Let  Aut$_0S$ be the identity component of Aut\,$S$, and 
Aut$_1S$ the normal subgroup of Aut\,$S$ of elements which leave each $C_i$ 
invariant.  By \cite{hau} we know that 
Aut$_0S \cong \C^*$, the multiplicative group of nonzero complex numbers. 
Our purpose is thus to determine the discrete part of Aut\,$S$. 
We shall summarize our results in Theorem \ref{mt}, Corollary \ref{mt} and 
the ensuing Remark below.   
In this note $\mu_m$ will denote the cyclic group of order $m$. 

\begin{Theorem}\label{mt}
Let  $S$  be a parabolic Inoue surface with second betti number $m > 0$.  
Then we have the following: 
\begin{enumerate} 
\item 
Aut$_0S$ $(\cong \C^*)$ coincides with the center of Aut\,$S$. 
\item 
Aut$_1S$ is commutative and is isomorphic 
to $\C^*\times \mu_m$.  
\item 
Aut\,$S$ is isomorphic to the semidirect product $\mu_m \ltimes \rm{Aut}_1S$, 
where the action of a generator of $\mu_m$ maps $(t,s)\in \rm{Aut}_1S  
\cong \C^*\times \mu_m$ to $(st,s)$, where $s\in \mu_m$ is considered as an $m$-th root 
of unity. 
\end{enumerate}  
\end{Theorem} 

From this we get easily the following: 

\begin{Corollary}\label{mc}  
The following hold: 
\begin{enumerate} 
\item 
There exist precisely $m$ cyclic groups of order $m$ in Aut$_1S$  
which have trivial intersection with Aut$_0S$.  
These are conjugate to one another in Aut\,$S$. 
\item 
Aut\,$S/\rm{Aut}_0S$ is isomorphic to the abelian group  $\mu_m\times \mu_m$. 
\item 
Aut\,$S/\rm{Aut}_1S$ is isomorphic to the cyclic group   $\mu_m$.  
There exist precisely $m$ cyclic groups of order $m$ in Aut\,$S$  
which is mapped isomorphically onto this quotient. 
\end{enumerate}  
\end{Corollary} 

{\em Remark}. 
\begin{enumerate} 
\item 
$S$  has the infinite cyclic fundamental group and 
its universal covering space  $W$  is a toric surface with a natural action of 
the algebraic torus $G:=\C^{*2}$.  
Then Aut$_1S$ is naturally identified with the maximal subgroup of $G$ 
whose action commutes with the covering transformations. 
Every element of  Aut\,$S$ can be explicitly described in terms of its lift 
to Aut\,$W$ with respect to 
the toric coordinates on  $W$ up to covering transformations. 
\item 
Let $D_m$ be the dihedral group of order $2m$, 
considered as the symmetry group of the graph $\Gamma$ of the cycle $C$, 
which a regular $m$-gon.  Then the action of Aut\,$S$ on $C$  induces 
a natural homomorphism  $u:$ Aut\,$S \ra  D_m$ with kernel Aut$_1S$. 
This is surjective onto the unique cyclic subgroup of order $m$. 
In particular, 
there exists no automorphism of $S$ which induces a reflection of $\Gamma$. 
\item 
Any of the $m$ cyclic groups of 3) of Corollary \ref{mc} acts on  $S$  without fixed 
points and permutes the irreducible components of  $C$  cyclically. 
Correspondingly, we have $m$ quotient surfaces $S_i, 1\leq i\leq m$, 
which are parabolic Inoue surfaces with second betti number  $1$.  
\item 
Since the homology classes of $E$ and $C_i$ generate $H_2(S,\Z)$ over  $\Q$  and 
$H_2(S,\Z)$ is torsion free,  Aut$_1S$ is precisely the subgroup of Aut\,$S$ of elements 
which act trivially on  $H_2(S,\Z)$. 
\item 
Parabolic Inoue surfaces with fixed second betti number $m$ are 
naturally parametrized by complex numbers $\al $ with $0<|\al |<1$ as $S=S(m,\al)$.  
The structure of the automorphism group 
of these surfaces are independent of the parameter  $\al$.  It follows that 
$\{S(m,\al)\}_\al$ actually gives the moduli space of parabolic Inoue surfaces.  
\item 
The Zariski open set 
$U:=S-C$ is invariant by the action of Aut\,$S$ as remarked above.  
$U$ has the natural structure of a holomorphic line bundle of degree $-m$ 
over the elliptic curve $E$ with $E\subseteq U$ the zero section (cf.\ \cite{e}).  
It turns out that Aut\,$S$ is isomorphic to the group of 
bundle automorphisms of  $U \ra E$  which induce translations on  $E$, 
the {\em theta group} in the sense of \cite{mu}.  The structure of the latter group 
is given in \cite{mu}, which gives in fact most of the structures of Aut\,$S$  
as stated in Theorem \ref{mt}  and Corollary \ref{mc}. 
\end{enumerate}  

\vspace{3 mm} 
The proof of the above results will be given in the next section. 
After some generalities on the behaviour of parabolic Inoue surface 
under Galois coverings in Section 3, 
we give in Section 4 some geometric description of the quotients of $S$ by 
typical subgroups of Aut\,$S$ in Theorem \ref{mt}.

Our interest on automorphism groups of parabolic Inoue surfaces 
comes from the question as to which 
anti-self-dual bihermitian structures on hyperbolic and parabolic Inoue surfaces 
constructed in \cite{fp} are invariant by the automorphisms of the surfaces (cf.\ 
also \cite{fjo}.

\section{Proof of Theorem \ref{mt}} 
A compact connected complex surface is said to be {\em of class} VII 
if its first betti number equals one and its Kodaira dimension equals $-\infty$.  
It is called {\em of class} VII$_0$ (resp.\ VII$_0^+$)  
if further it is minimal (resp.\ minimal and with positive second betti number).  

A parabolic Inoue surface is a surface of class VII$_0^+$ discovered by 
Inoue in \cite{ino} (cf.\ also \cite{na84}).  
A parabolic Inoue surface $S$ is written uniquely 
as  
\[ S=S(m,\al), 0<|\al |<1,\]  
where $m$ is the second betti number of  $S$ (cf.\ \cite[Lemma 3.6]{fp}). 
We explain its structure according to \cite{na84} (cf.\ also \cite{e}).  
(Note however the slight deviation of 
the notations from those of \cite{na84}.) 
$S$  has an infinite cyclic fundamental group, and 
its universal covering space $W$ is common to all  $S$.  

$W$ is covered by coordinate open subsets 
$U_k=\C^2(x_k,y_k),k\in \Z$, and  $V=\C^*(w)\times \C(x)$ 
with the transition relations  
\begin{gather}\label{tr}
 x_{k+1}=y^{-1}_k,\ \  y_{k+1}=x_ky_k^2 \ \  \mbox{on} \ U_k\cap U_{k+1}\\ 
w=x_ky_k,\ \  x=x_k^{k+1}y_k^k, \ \ \mbox{or} \ \ x_k=w^{-k}x,\ \  y_k=w^{k+1}x^{-1} 
\ \  \mbox{on}\ U_k\cap V \label{tr2}
\end{gather} 
$W$ is a toric surface on which the algebraic two-torus $G:=\C^*(s)\times\C^*(t)$ acts by 
\begin{gather}\label{ac}
(x_k,y_k) \ra  (s^{-k}tx_k,s^{k+1}t^{-1}y_k) \ \ \mbox{on} \ U_k \\ \mbox{and} \\   
(w,x) \ra  (sw,tx)\ \  \mbox{on} \ V. \label{ac2} 
\end{gather} 
The $\C^*(t)$ action fixes pointwise the curve  $\tilde{E}:=\{x=0\}\cong \C^*$ in $V$. 
The equations  $x_k=0$ on  $U_k$ and $y_{k+1}=0$ on  $U_{k+1}$ define 
a compact smooth rational curve $\tilde{C}_{k+1}$ on  $W$: 
\begin{equation}\label{def} 
\tilde{C}_{k+1}: \ 
x_k=0\ \ \mbox{on} \ U_k, \ \  y_{k+1}=0\ \ \mbox{on}\  U_{k+1}.  
\end{equation}  

Now for any complex number $\al, 0<|\al |<1$, and for any positive integer  $m$ 
we define an automorphism $\gamma = \gamma_{m,\al}$ of $W$ which 
maps $U_{k-m}$ to $U_k$ and preserves  $V$ by 
\begin{gather}\label{bi}
(x_{k-m},y_{k-m}) \ra  (x_k,y_k)=(\al ^{-k}x_{k-m},\al ^{k+1}y_{k-m})\ \ \mbox{on}\ U_{k-m} \\ 
(w,x) \ra  (\al w,w^mx)  \ \ \mbox{on}\ V. \label{bi2} 
\end{gather} 
Note that for $m=0$ (\ref{bi}) and (\ref{bi2}) coincide with 
the $\C^*(s)$-action given in (\ref{ac}) and (\ref{ac2}) for $s=\al$. 
Then  $\gmm_{m,\al}$ generates a properly discontinuous, cocompact and fixed-point free 
infinite cyclic group of transformations of $W$.  
The quotient $W/\langle \gamma_{m,\al} \rangle $ 
is by definition a parabolic Inoue surface denoted by  
		\[S=S(m,\al)  \]
with second betti number $m$. 
The universal covering of $S$ is  $W$,  independently of $\al$ and $m$. 
By (\ref{bi2}) and (\ref{ac2}) $\gamma_{m,\al}$ normalizes $G$ 
in Aut\,$W$ by the formula  
\begin{equation}\label{sdr} 
 \gamma_{m,\al} (s,t)\gamma_{m,\al} ^{-1} = (s,s^mt). 
\end{equation}  
In particular $\gmm_{m,\al}$ commutes with $\C^*(t)$-action, 
and for the $\C^*(s)$-action we have 
\begin{equation}\label{sdq} 
s^{-1}\gamma_{m,\al}s = (1,s^m)\gamma_{m,\al}.  
\end{equation}  
Thus an element $(s,t)\in G$  commutes with $\gamma_{m,\al}$ if and only if 
$s$  is an $m$-th root of unity. 
Fix an $m$-th root $\beta $ of $\al$. 
We set $\gmm_\bt := \gmm_{1,\bt}$. Then by (\ref{sdq})  applied to the case $m=1$ 
and $\al = \bt$  
we get: 
\[  s\gmm^m_\bt s^{-1} = (1,s^{-m})\gamma_{\bt}^m.  \]  
Since the action of $\gmm^m_\bt$ on  $V$ takes the form 
\begin{equation}\label{itr} 
 (w,x) \ra  (\bt^mw, \bt^{m(m-1)/2}w^mx),  
\end{equation}  
if we take  $s=\bt^{(m-1)/2}$ (with any one of the two values fixed when $m$ is even), 
by (\ref{bi2}) we obtain 
\begin{equation}\label{sdn} 
 \nu^m = s\gmm^m_\bt s^{-1} = \gamma_{m,\al}, \ \ 
\end{equation}  
where 
\begin{equation}\label{ndef} 
\nu=\nu_{\bt} := s\gmm_\bt s^{-1}
\end{equation}  
which acts in the following form on $V$:  
\[\nu : \ (w,x)\ra (\bt w, \bt^{-(m-1)/2}wx). \ \ \]
\indent  
We summarize the implications of the above computations in the following: 
\begin{Lemma}\label{sum}
Aut\,$S$ contains a subgroup $H$ isomorphic to a semidirect product 
$\mu_m \ltimes (\mu_m\times \C^*)$ as in 3) of Theorem \ref{mt} realized 
in the form  $H=\langle \bar{\nu} \rangle \ltimes 
(\langle \rho \rangle \times \C^*(t))$, where $\bar{\nu}$ is the 
automorphism of $S$  induced by $\nu$, and 
$(\rho, t)\in G, \rho =\exp(2\pi \sqrt{-1}/m)$, admits a naturally induced 
action on  $S$. 
\end{Lemma} 

{\em Proof}.  By (\ref{sdn}) 
$\nu$  induces on $S$  a fixed point free automorphism $\bar{\nu}$ 
of order  $m$. 
If $C_i$ is the natural image of the curve $\tilde{C}_i,1\leq i\leq m$, in $S$, 
$C=C_1+\cdots +C_m$ is the unique cycle of rational curves on $S$. 
Since the action of $\nu$  maps $\tilde{C}_k$ to $\tilde{C}_{k+1}$, 
$\bar{\nu}$  transforms cyclically the curves $C_i$. 
Thus we have $\langle \bar{\nu} \rangle \cap$ \rm{Aut}$_1S = \{e\}$, 
where  $\langle \bar{\nu} \rangle $ is the cyclic group 
of order $m$  generated by $\bar{\nu}$ and $e$  is the identity of  Aut\,$S$. 
On the other hand, by (\ref{sdr}) 
the action of $(\rho ,t) \in  G$  commutes with  $\gmm_{m,\al }$; hence 
we may consider $\langle \rho \rangle \times \C^*(t)$ as a subgroup of  Aut\,$S$, 
which clearly is contained in Aut$_1S$.  In view of (\ref{sdr}) for $m=1$ 
these two groups form 
a semi-direct product $H:= \langle \bar{\nu} \rangle \ltimes 
(\langle \rho \rangle \times \C^*(t))$ in Aut\,$S$ as in 3) of Theorem \ref{mt}.  
The identity component of $H$ is clearly isomorphic to $\C^*(t)$ and  
the quotient group $H/\C^*(t)$ is isomorphic to $\mu_m\times \mu_m$, again 
by (\ref{sdr}). \hfill$\square$ 

\vspace{3 mm} 
It remains to show that $H$ coincides with Aut\,$S$. 
First we prove the following: 

\begin{Lemma}\label{h1}
Let  $h$  be any automorphism of finite order on  $S$.  
Then the restriction $h_E$ of  $h$  to  $E$  is a translation. 
\end{Lemma} 

{\em Proof}.  
First we note that the induced action of  $h$  on  $H^1(S,O_S)\cong \C$ is 
trivial.  Indeed, let  $\hat{S}$ be the quotient of $S$ by the cyclic group $H$ generated 
by  $h$  and $\tilde{S}$ be a resolution of $\hat{S}$.   
We have the natural isomorphisms  
$H^1(S,O_S)^H\cong H^1(\hat{S},O_{\hat{S}})\cong H^1(\tilde{S},O_{\tilde{S}})$, 
where $(\ )^H$ denotes the subspace of $H$-fixed elements. 
If the action on  $h$  on $H^1(S,O_S)$ is non-trivial, 
$H^1(\tilde{S},O_{\tilde{S}})\cong H^1(S,O_S)^H = \{0\}$. 
Then by Kodaira's classification of surfaces,  $\tilde{S}$  must be a K\"ahler surface. 
But then  $S$ also is K\"ahler, as it is bimeromorphic to a finite covering 
of a K\"ahler surface  $\tilde{S}$.  This is a contradiction.  Thus  $h$  acts trivially 
on  $H^1(S,O_S)$.  Since 
the restriction map  $H^1(S,O_S) \ra H^1(E,O_E)$ is known to be isomorphic 
and $h$-equivariant, this implies that the action of $h_E$ on $H^1(E,O_E)$ also is trivial. 
This implies that  $h_E$ is a translation. 
\hfill $\square$

\vspace{3 mm} 
In passing 
we ask the following question suggested by the above lemma. 
Let  $S$  be any compact complex surface in class VII.  
Does any automorphism $h$ of  $S$  act trivially on  $H^1(S,O_S)\cong \C$ ? 
By the same proof as above this is true when  $h$  is of finite order. 
Does  $S$  admit any automorphism of infinite order if $\dim$ Aut\,$S=0$? 

\vspace{3 mm} 
Now 
$U:=S-C=V/\langle \gmm_{m,\al} \rangle $  has a natural structure of a holomorphic line bundle 
$u: U \ra  E$ over the elliptic curve  $E$  such that 
$E\subseteq U$  is the zero section of $u$ and 
the induced $\C^*$-action on  $U$  is the natural fiber action 
on this line bundle.  
Now let  $g$   be any automorphism of  $S$, which necessarily preserves  $U$ and 
$g$  normalizes $\C^*=$Aut$_0S$.  
\begin{Lemma}\label{cnt}
The above $g$ actually  centralizes  $\C^*$ 
so that $g$  acts on  $U$  as a bundle automorphism of  $u$.   
\end{Lemma} 

{\em Proof}.  
Suppose that  $g$  does not centralize  $\C^*$. 
Then for any $t\in \C^*$ and $p\in S$ we have $g^{-1}tg(p) = t^{-1}(p)$. 
Take  $p$  with $p\in U-E$ and let $t$ tend to zero. 
Then the left hand side tends to a point on $E$, 
while the right hand side has no limit in $U$ since  $t^{-1}(p)$ goes 
to infinity on the fiber over  $u(p)$, which is a contradiction.  
\hfill $\square$ 

\vspace{3 mm} 
Since the degree of $U$ equals $-m=E^2<0$, 
the group $H(U)$ of bundle automorphisms of $u$ is 
an extension by $\C^*$ of a finite group of automorphisms of  $E$ 
(cf.\ Lemma \ref{theta} below). 
\begin{Lemma}\label{bat}
Let  $g$  be as above. 
The induced automorphism  $g_E$  of $E$  is a translation.  
\end{Lemma} 

{\em Proof}.  
If  $g$  is of finite order, this follows from Lemma \ref{h1}.  Thus 
it suffices to show  that the composition $tg$ becomes 
of finite order for some $t\in\C^*=$\rm{Aut}$_0S$ since  $(tg)_E = g_E$.  
Indeed, since $g_E^k = 1_E$, the identity of $E$, for some $k>0$, and $u$  is equivariant, 
$g^k$  is a bundle automorphism of $U$  over the identity of $E$, and 
hence  $g^k$ belongs to $\C^*$.  
Take any element $t$ of  $\C^*$  with $t^k=(g^k)^{-1}$. 
Then  $tg$  clearly is of order $k$ as desired. 
\hfill $\square$ 

\vspace{3 mm} 
Let  $G(U)$ be the subgroup of $H(U)$ of 
automorphisms of $U$  which induce 
on $E$ a translation.  We summarize the known structure of $G(U)$ from \cite{mu} 
in the following: 

\begin{Lemma}\label{theta}
$G(U)$ fits into the exact sequence 
\begin{equation}\label{th} 
  0 \ra  \C^* \ra  G(U) \ra  K(U) \ra  0
\end{equation}  
of algebraic groups 
in the sense that $G(U)$ is a central extension of $K(U)$ by $\C^*$, 
where $K(U)$ is the finite subgroup of the group of translations $b$ of $E$ 
consisting of those $b$ with  $b^*U\cong U$, 
while $\C^*$ corresponds to the natural fiberwise $\C^*$-action on $U$.  
Moreover, $K(U)$ is isomorphic to  $\mu_m\times\mu$, and 
$\C^*$ coincides with the center of  $G(U)$. 
\end{Lemma} 

{\em Proof}.  
These are all found in \cite{mu} 
where in general line bundles over an abelian variety of arbitrary dimension are treated.  
Specilized to one-dimensional case we easily deduce the results stated above: 
See \cite[p.225,Th.1]{mu} for the sequence (\ref{th}); \cite[p.60,Appl.1]{mu} 
for the finiteness of $K(U)$; \cite[p.84,(v)]{mu}, \cite[p.154 Remarks]{mu} and 
\cite[p.150, R-R, p.155 (*)]{mu}) for its computation; and finally, 
\cite[p.155,Cor.2]{mu} (cf.\ \cite[p.223,Cor.]{mu}) for the last assertion.   
\hfill $\square$ 

\vspace{3 mm} 
{\em Proof of Theorem \ref{mt}}.  
So far we have shown that  
$H:= \langle \bar{\nu }_\beta \rangle \ltimes 
(\langle \rho \rangle \times \C^*(t))\subseteq {\rm Aut} S \subseteq K(U)$ 
and  $(\langle \rho \rangle \times \C^*(t))\subseteq {\rm Aut}_1S$, the second 
inclusion being due to Lemmas \ref{cnt} and \ref{bat}. 
From the first inclusion we get the natural inclusion $H/\C^*(t)\subseteq G(U)$,  
while both groups are isomorphic to  $\mu_m\times \mu_m$ by Lemma \ref{theta}. 
Thus we have the equalities  $H={\rm Aut} S = G(U)$.  
All the assertions of Theorem \ref{mt} then follows from the structure of 
$H$  already mentioned or of the structure of $G(U)$ given in Lemma \ref{theta}. 
\hfill $\square$ 

\vspace{3 mm} 
The statements in Remark after Corollary \ref{mc} are immediate from 
our arguments above.  
Instead of using the structure of the theta group as in the above proof 
one can also obtain the same results by direct computations, which however we shall 
omit here.

\section{Galois Coverings and parabolic Inoue surfaces} 

It seems that parabolic Inoue surfaces are closed under finite coverings.  
Namely we could ask if the following is true: 
Let  $S_1 \ra S_2$ be 
a generically surjective meromorphic map 
of compact complex surfaces of class VII$_0^+$. 
Then $S_1$ is a parabolic Inoue surface if and only if  so is $S_2$.  

In this section, however, we consider just the simplest cases where  $f$  is 
a (holomorphic) Galois covering.  For this purpose the following 
characterization of a parabolic Inoue surface due to 
Hausen \cite{hau} is useful. 

\begin{Lemma}\label{hau}
Let  $S$  be a surface of class VII$^+_0$.  Then 
$S$  is a parabolic Inoue surface if and only if 
$S$ admits an effective $\C^*$-action.  
\end{Lemma} 

One could have also used the characterization of by \cite{na84}: 
a parabolic Inoue surface is precisely a surface of class VII$^+_0$ 
which carries a smooth elliptic curve and a cycle of rational curves.  But 
the argument seems easier with the above characterization. 

\begin{Lemma}\label{chp} 
Let  $S$  be a parabolic Inoue surface and  $G$  a finite subgroup of Aut\,$S$.  
Let  $S':=S/G$  be the quotient surface, and 
$\bar{S}$ the (unique smooth) minimal model which is bimeromorphic to 
$S'$.   Then  $\bar{S}$  is again a parabolic Inoue surface. 
\end{Lemma} 

{\em Proof}.  Since  $G$  centralize  Aut$_0S\cong\C^*$ by Theorem \ref{mt}, 
the $\C^*$-action on  $S$  descends to $S'$, 
then lifts to its minimal resolution $\bar{S}'$, and then descends to $\bar{S}$.  
On the other hand, 
since  $G$  preserve the unique elliptic curve $E$  on  $S$  and acts by translations 
on  $E$,  $S'$, and then also $\bar{S}'$ and  $\bar{S}$, contain a nonsingular elliptic curve 
on it.  Moreover, if  $\bar{E}$ is the elliptic curve on  $\bar{S}$, 
its self-intersection number is negative as well as that of $E$.  
Hence the second betti number of $\bar{S}$  
is positive and $S$ belongs to class VII$^+_0$. 
Therefore by Lemma \ref{hau} $\bar{S}$  is a parabolic Inoue surface.  
\hfill $\square$ 

\vspace{3 mm} 
The next lemma gives also the propagation of parabolic Inoue property  
but in the converse direction. 

\begin{Lemma}\label{rcv}
Let  $\bar{S}$  be a parabolic Inoue surface.  
Let  $u: \hat{S}\ra \bar{S}$  be a finite Galois covering with Galois group $H$ 
where $\hat{S}$ is a normal complex surface.  
Let  
$S$ be the (unique smooth) minimal model which is bimeromorphic to 
of $\hat{S}$.  
Then  $S$ is again a parabolic Inoue surface. 
\end{Lemma} 

{\em Proof}.  
Let  $\bar{C}$ be the discriminant locus of  $u$  on  $\bar{S}$ and 
$\hat{C}=u^{-1}(\bar{C})$ with reduced structure. 
Since  $\bar{C}$ has only normal crossings, 
$\hat{S}$  has only cyclic quotient singularities.  Thus we have 
the naturally defined notion of 
the sheaf $\Theta_{\hat{S}}(-\log \hat{C})$ 
of logarithmic tangent vector fields on  $\hat{S}$  along $\hat{C}$ 
(cf.\ \cite{st}). 
On the smooth part $\hat{S}_0$ of  $\hat{S}$  the natural sheaf homomorphism 
$\Theta_{\hat{S}}(-\log \hat{C}) \ra u^*\Theta_{\bar{S}}(-\log \bar{C})$ is $H$-equivariantly 
isomorphic and hence the subspaces of $H$-fixed elements are also isomorphic. 
Thus any holomorphic vector field which comes from the $\C^*$ action 
gives rise to a section of $u^*\Theta_{\bar{S}}(-\log \bar{C})$ on $\hat{S}_0$; 
then considered as a section of $\Theta_{\hat{S}}(-\log \hat{C})$ 
by the above isomorphism, 
it extends by the normality of $\hat{S}$ 
to a section of $\Theta_{\hat{S}}(-\log \hat{C})$ on the whole  $\hat{S}$.  
Then it lifts to its minimal 
resolution $S'$ and then descends to the minimal model  $S$.  Thus  $S$ admits 
the induced $\C^*$-action.  Moreover, it is easy to show that the proper transform 
of the inverse image of the unique elliptic curve $\bar{E}$  on  $\bar{S}$  has 
a negative self-intersection on  $S$, and hence that $S$ is in class VII$_0^+$. 
Thus by Lemma \ref{hau} it is a parabolic Inoue surface. \hfill $\square$

\section{Typical cyclic subgroups}	

In this section we shall identify the minimal model 
of the quotients of a parabolic Inoue surface $S$ 
by the typical finite cyclic subgroups of Aut\,$S$ given in Theorem \ref{mt}.  
In deciding the isomorphism classes of the above minimal modeles, 
the following lemma is useful.  

\begin{Lemma}\label{36}
Let  $S$  be a parabolic Inoue surface.  Let  $u: W \ra S$  be the universal covering 
map of  $S$  and  $\tilde{E}=u^{-1}(E)$, where $E$ is the unique elliptic curve on $S$.  
Suppose that $E^2=-m$ and $\tilde{E} \ra  E$  is isomorphic to the natural quotient map 
$\C^* \ra  \C^*/\langle \al\rangle $, then  $S$ is isomorphic to  $S(m,\al)$.  
\end{Lemma} 

{\em Proof}.  The $m$-part is clear from the description of Section one.  
For the $\al$-part see Lemma 3.6 of \cite{fp}. \hfill $\square$ 

\vspace{3 mm} 
In Examples below we assume that $S=S(m,\al), m\geq 1, 0<|\al|<1$.  

\vspace{3 mm} 
\noindent{\bf Example 1}. 
The $m$ cyclic groups $H_\bt$ of order $m$ with trivial intersection 
with  Aut$_1S$ mentioned in 3) of Corollary \ref{mc} are 
parametrized by the set  $B=B(\al)$ of $m$-th roots of $\al$ as 
\[
H_\bt:= \langle \bar{\nu}_{\bt}\rangle, \ \  \bt \in B 
\] 
with $\nu_{\bt}$ as in (\ref{ndef}). 
The quotient $S_\bt:= S/H_\bt$ 
fits in an unramified covering  $S \ra S_{\bt} \cong S(1,\bt)$.  
Thus in all we have $m$ such unramified quotients of  $S$. 

\vspace{3 mm} 
\noindent{\bf Example 2}. 
Let  $\triangle_l$  be 
the unique cyclic subgroup of order $m$ in $\C^*(t)=$ Aut$_0S$.  
Let $S':= S/\triangle_l$ be the quotient  of  $S$  by $\triangle_l$  and 
$\bar{S}$ the smooth surface obtained by taking the minimal resolution of  
the singularities of $S'$.  Then  $\bar{S}$ is a parabolic Inoue surface 
isomorphic to  $S(ml,\al)$.  

\vspace{1 mm} 
{\em Proof}.  
The action of $\triangle_l=\langle\rho  \rangle, \rho = \exp(2\pi\sqrt{-1}/l)$, 
fixes  $E$ pointwise and its 
action on $U_k(\subseteq W)$ takes the form 
\begin{equation}\label{say} 
 (x_k,y_k) \ra  (\rho x_k,\rho ^{-1}y_k). 
\end{equation}  
Thus the quotient  $S'$ 
has $m$ rational double points $p'_i$ of type A$_{l-1}$, namely the images of 
the nodes  $p_i:=C_i\cap C_{i+1}, 0\leq i\leq m-1$, of $C$, 
and is smooth otherwise. 
The minimal resolutions at each $p'_i$ give rise to a chain of $(l-1)$ $(-2)$-curves.  
Together with the proper transforms $\bar{C}'_i$ in  $\bar{S}$ 
of images of irreducible components  $C_i$  of $C$ in $S'$ these can be shown to 
form a cycle $\bar{C}$ of rational curves consisting of  $m(l-1)+m =ml$ irreducible components. 
Other than these $ml$ curves 
the only curve on the resoloved surface  $\bar{S}$ is 
the inverse image $\bar{E}$ in $\bar{S}$ of the image of $E$ in $S'$.  
Since $\triangle_l$ pointwise fixes the elliptic curve $E$, we conclude 
easily that $\bar{E}\cong E$ and 
that the self-intersection number of $\bar{E}$ equals $-lm$.  

Now if $\bar{S}$ is not minimal, we obtain its minimal model $T$ contracting 
some of the irreducible components of $\bar{C}$ to points.  By 
Lemma \ref{chp} $T$ is a parabolic Inoue surface.  This, however, implies that 
the image of $\bar{C}$ in  $T$  must be a cycle of rational curves with 
precisely $ml$ irreducible components, since the contraction is isomorphic 
at the neighborhood of  $\bar{E}$.  This is impossible as  $\bar{C}$ contains 
only $ml$ irreducible components. 
Hence $\bar{S}$ itself is minimal and 
is a parabolic Inoue surface. 
Since it contains an elliptic curve 
with self-intersection number  $-lm$, it is of the form  $S(ml,\al')$ for 
some $\al'$. 
 
 It remains to show that $\al '=\al $.  For this purpose we apply to  $W$ 
the construction of $S'$ and $\bar{S}$  starting from  $S$.  
Namely we take the quotient $W'=W/\triangle_l$  and take its minimal resolution 
$\bar{W}$.  Since the action of  $\C^*$ commutes with 
$\gamma = \gamma_{m,\al}$, the action of $\gamma$ descends to  $W'$  and then 
lifts to  $\bar{W}$.   
The resulting action of  $\bar{W}$  is again properly discontinuous and fixed point 
free and the quotient  $\bar{W}/\langle \gamma\rangle$  is isomorphic to  $\bar{S}$. 
In particular,  $\bar{W}$  is isomorphic to $W$.   
By construction, the inverse image 
of $\bar{E}$  in  $\bar{W}$  is naturally and $\gamma$-equivariantly identified 
with  $\{x=0\}\cong \C^*(w)$ in (the original) $W$.   
Since the action of $\gmm$ on  $\C^*(w)$ 
is by multiplication by  $\al $, the assertion follows by Lemma \ref{36}. 
\hfill $\square$

\vspace{3 mm} 
\noindent{\bf Example 3}. 
Recall that $S=S(m,\al)$. 
The $m$ cyclic groups in  Aut$_1S$ which have trivial intersections in 
Aut$_0S$ are identified with the subgroups, say $M_j$, of $G=\C^{*2}$ 
generated by  $g_j:=(\rho ,\rho ^j), 0\leq j\leq m-1, \rho = \exp(2\pi \sqrt{-1}/m)$.  
The actions of $g_j$ on $U_k$ and $V$ take respectively the forms:  
\begin{equation}\label{bb} 
g_j: (x_k,y_k) \ra  (\rho ^{-k+j}x_k,\rho ^{k-j+1}y_k), \ \ (w,x) \ra (\rho w,\rho^jx). 
\end{equation}  
Let  $H$  be any of these subgroups 
and   $l$  any integer $> 1$ which divides  $m$. 
Let  $H_l$ be the unique cyclic subgroup of  $H$ of order $l$.  
Since the above  $m$ groups are conjugate to each other in Aut\,$S$, 
the quotient  $S \ra S':=S/H_l$  is up to isomorphisms independent of  $H$.  
By (\ref{bb}) the action of $H_l$ on  $S$  
is fixed point free outside the curve  $C$; 
the points  $p_k:=C_k\cap C_{k+1}$ are all fixed, but for instance 
the irreducible components $C_{li+1},\ 0\leq i < n, n:=m/l$,  
are also pointwise fixed.  
So the singularities of the quotient  $S/H_l$ is not easily identified. 
In what follows, 
by describing this quotient in another way 
we identify the (unique smooth) minimal model  $\bar{S}$  which is bimeromorphic 
to  $S'$.  Namely we show the following: 
\begin{Proposition}\label{al}
We have $\bar{S} \cong S(n,\al^l)$. 
Moreover, the natural bimeromorphic map 
$S'\ra \bar{S}$ is holomorphic and obtained by contracting the image of 
each connected component of $C-\bigcup_{0\leq i < n}C_{li+1}$ to a point. 
\end{Proposition} 
For the proof 
we start in general a parabolic Inoue surface  $\bar{S}=S(n,\beta)$ 
with a cycle of rational curves $\bar{C}=\bar{C}_1+\cdots +\bar{C}_n$ and 
with a smooth elliptic curve $\bar{E}$ with $\bar{E}^2=-n$. 
The line bundle $[\bar{C}]$ defined by  $\bar{C}$ belongs to  
Pic$_0\bar{S}\cong \C^*$ corresponding to the number $\beta\in \C^*$ 
(cf.\ \cite[p.425]{na84}).  
Indeed, by (\ref{def}) the defining equation $f=0$ 
of $\tilde{C}:=\bigcup_k\tilde{C}_k$  on  $W$  is given by $x_ky_k=0$ on  $U_k$ and 
$w=0$ on $V$, and by (\ref{bi2}) 
$\nu_\bt$ acts on  $f$  by $\nu_{\bt}^*f=\beta f$.  
This shows the assertion. 

Hence 
$[\bar{C}]$ is divisible by $l$ in  Pic$_0\bar{S}$ (exactly $l$ roots) and 
thus we obtain $l$ ramified cyclic coverings of degree $l$ 
which are totally ramified along  $\bar{C}$ and unramified otherwise.  
Let $w: \hat{S} \ra  \bar{S}$  be one of them and  $L$  the corresponding holomorphic line bundle 
with $L^l = [\bar{C}]$. Then $L^j,1\leq j< l$, are never trivial when restricted to the elliptic curve 
 $\bar{E}$; indeed, the restriction map Pic$_0\bar{S} \ra$ Pic$_0\bar{E}$ is 
identified with the quotient $\C^* \ra \bar{E}$  with 
its kernel generated by  $\beta$ and hence by  $[\bar{C}]$.  
This implies in particular 
that $\hat{E}:=w^{-1}(\bar{E})$ is connected and with self-intersection 
number $-ln$. 

$\hat{S}$ has $n$ rational double points of type A$_{l-1}$, and 
by taking the minimal resolutions of these singularities 
we obtain from each of them $(l-1)$ rational curves 
with self-intersection number $-2$.  Let  $v: S\ra \hat{S}$ be the 
minimal resolution. 
For $1\leq i\leq n$ 
let $\hat{C}_i$ be the inverse images of $\bar{C}_i$ in  $\hat{S}$ 
(with reduced structure),
and  $C_i$ their proper transforms in $S$.  Then we can show that 
altogether 
we get a cycle of $n+n(l-1)=nl$ rational curves on $S$.  
Now we conclude that 
$S$ is again a parabolic Inoue surface with second betti number $nl$ as follows.  

\begin{Lemma}\label{prb}
Let $\bar{S}=S(n,\bt)$ and $S$ be as above.  Then 
$S$ is a parabolic Inoue surface of the form  $S(nl,\al)$ for some $\al$ with 
$\al^l=\bt$. 
\end{Lemma} 

{\em Proof}.  
We first show that 
$S$ is a parabolic Inoue surface of the form  $S(nl,\al)$ for some $\al$.  

By Lemma \ref{rcv} 
the minimal model $S'$ of  $S$  is 
a parabolic Inoue surface. 
Since the blowing-down map $u: S\ra S'$ is isomorphic in a neighborhood of  
the proper transform  $E$  of  $\hat{E}$, the self-intersection number of 
the image of $E$ in $S'$ is again equal to $-nl$.  Thus  $S'$ should contain 
a cycle of $nl$ rational curves. But 
by construction  $S$ contains no curves other than $E$  
and the $nl$ rational curves mentioned above.  Thus no curve on  $S$  cannot be 
blowing down to a point of  $S'$; namely  $S$  already is minimal and is a parabolic 
Inoue surface of the form  $S(nl,\al)$.  

It remains to show that $\al^l=\bt$.  
We first construct the universal covering $b: \tilde{W} \ra  S$ of $S$.  
Let $a: \bar{W} \ra  \bar{S}$  be the universal covering of  $\bar{S}=S(n,\bt)$ 
so that $a^{-1}(\bar{E}) \ra  \bar{E}$ is isomorphic 
to the quotient $\C^*\ra \C^*/\langle \bt\rangle $.  
Then the pull-back $b: \tilde{W}:=\bar{W}\times_{\bar{S}}S \ra  S$ of $a$  
to  $S$ {\em via}  $wu: S \ra  \bar{S}$ is an infinite cyclic unramified covering of  $S$.  
Since  $S$  is a parabolic Inoue surface and 
hence has an infinite cyclic fundamental group,  
$b$  must be the universal covering of  $S$ and   
$b^{-1}(E)$ is isomorphic to  $\C^*$.  
Moreover, the projection  $c: \tilde{W} \ra  \bar{W}$  induces 
an unramified cyclic covering $c': b^{-1}(E) \ra a^{-1}(\bar{E})$ of degree $l$.  
Thus $c'$ is isomorphic to the map  $\C^* \ra  \C^*, s \ra  s^l$,  and 
the multiplication by $\al$ is sent to that by $\al^l$ on the image $\C^*$.  
Thus by Lemma \ref{36} we have  $\al ^l=\bt$. 
\hfill $\square$ 

\vspace{3 mm} 
{\em Proof of Proposition \ref{al}}.  
By Lemma \ref{prb} we may assume that  $S$  and  $\bar{S}$  are the surfaces in that lemma 
with $m=ln$. Indeed,  $n$ and $\bt$ in $\bar{S}(n,\bt)$ of the lemma 
can be chosen arbitrarily and the same is true for $\al$ with 
$\al^l = \bt$ by a suitable choice of  $\hat{S}$. 
Let  $S \stackrel{q}{\ra} S' \stackrel{r}{\ra} \bar{S}$ be 
the Stein factorization of  $wv: S \ra  \bar{S}$ so that 
we have the commutative diagram 
\begin{equation}\label{cmd} 
\begin{CD}
	  S	 @>v>>  \hat{S} \\
		  @VqVV	  @VVwV	      \\
	S'  	@>r>>  \bar{S} , 
\end{CD}                      
\end{equation}       
where $v$ and $w$ are the natural maps defined above. 
The action of the Galois group $K$ of $w: \hat{S} \ra  \bar{S}$ lifts to 
$S$  and  $q$  is identified with the quotient map $S \ra S/K$. 
The action of  $K$  preserves each irreducible components of the curves on 
$S$  and therefore 
$K$ is contained in  Aut$_1S$.  Moreover, since  $K$  acts freely 
on $E(\cong\hat{E})$  in  $S$, the intersection of $K$ and Aut$_0S$ 
is trivial.   Thus  $K$  must be contained in one of the maximal cyclic subgroups 
$M_j$.  Thus the first assertion is proved.  The second one follows immediately 
from the construction of $S$  from  $\bar{S}$.    
\hfill $\square$

\vspace{3 mm} 
\noindent {\bf Example 4}: We consider the special case of Example 3 for $l=2$.  
In this case we can describe the quotient directly. 
Let  $S=S(m,\al)$ as before with $m=2n$ even. 
Then we obtain an involution $\iota:=g_l^n$ (cf.\ (\ref{bb})). 
Explictly, this is induced by the involution  $\tilde{\iota }$  on  $W$ defined by 
\begin{gather}\label{iv}
\tilde{\iota }: \ \ (x_k,y_k) \ra ((-1)^kx_k,(-1)^{k+1}y_k)\ \ \mbox{on} \ U_k,\ \ 
 (w,x) \ra  (-w,x)\ \ \mbox{on} \ V. 
\end{gather} 
We have a cycle of rational curves $C=C_1+\cdots +C_{n}$ on $S$, 
where $C_i$ is the natural image of $\tilde{C}_i$.   
Then  $C_1+C_3+\cdots +C_{n-1}$ is precisely the fixed point set of $\iota$. 

Therefore 
the quotient  $S':=S/\langle \iota \rangle $ is a smooth surface, 
whose structure is described as follows.  
Let  $C'_i$ be the images of  $C_i$. 
Then we have  $(C'_{2k-1})^2=-4$ and  $(C'_{2k})^2=-1$. 
Thus we may contract $C'_{2k}$ to points 
and obtain another smooth surface $\bar{S}$.  
Let  $\bar{C}_i$ be the images of  $C'_i$ and $\bar{C}$ that 
of  $C'=\sum_i C'_i$ in $\bar{S}$.  
Then $\bar{C}$ is 
of the form  $\bar{C}=\bar{C}_1+\bar{C}_3+\cdots +\bar{C}_{2l-1}$ 
with $(\bar{C}_{2k-1})^2=-2$.  
The image $\bar{E}=E/\langle \iota \rangle $ of the elliptic curve $E$ on $\bar{S}$ 
is again a smooth elliptic curve. 
Clearly  $\bar{E}$ and  $\bar{C}$ are the unique curves on  $\bar{S}$, and 
$\bar{S}$ is again a parabolic Inoue surface and 
with second betti number $n$.  Thus the quotient is isomorphic to 
$S(n,\al^2)$. 

Let  $S \stackrel{v}{\ra}  \hat{S} \stackrel{w}{\ra} \bar{S}$ be 
the Stein factorization of  $S \ra  \bar{S}$.  
It turns out that $v$ is nothing but 
the contraction of the $l$ (-2)-curves $C_{2k}$ to $l$ ordinary double points 
on the normal surface $\hat{S}$, and 
$w$ is the branched double covering 
with branch locus $\bar{C}$. 
We thus recover the diagram (\ref{cmd}) from the other direction. 
\hfill$\square$

\vspace{3 mm} 
\begin{flushright}
Department of Mathematics\\
Graduate School of Science\\
Osaka University\\
Toyonaka 560-0043 Japan \\ 
{\small E-mail address: fujiki@math.sci.osaka-u.ac.jp} \\ 
\end{flushright}


{\footnotesize \begin{thebibliography}{99} 

\bibitem{e}
Enoki, I., 
 Surfaces of class ${\rm VII}_0$ with curves,
T\^ohoku Math. J., {\bf 33} (1981), 453--492.
\bibitem{fjo} 
Fujiki, A., 
Bihermitian anti-self-dual structures on compact non-K\"ahler surfaces, 
In: Proc. of Complex Geometry in Osaka, 
Lecture Notes in Math., {\bf 9}, 2008, 137-151.  
\bibitem{fp} 
Fujiki, A., and Pontecorvo, M., 
Anti-self-dual bihermitian structures on Inoue surfaces, 
arXive:math.DG/0903.1320v1, (2009). 
\bibitem{hau}
Hausen, J.,
Zur Klassifikation of glatter kompakter $\C^*$-Fl\"achen, 
Math. Ann. {\bf 301} (1995), 763--769.
\bibitem{ino} 
Inoue, M., 
New surfaces with no meromorphic functions,
Proc. of ICM (Vancouver, B. C., 1974), 423--426. 
\bibitem{mu} 
Mumford, D., 
Abelian varieties, 
Tata Institute and Oxford Univ. Press, 1974. 
\bibitem{na84}
Nakamura, I.,
 On surfaces of class ${\rm VII}_0$ with curves,
Invent. Math. {\bf 78}, (1984), 393--443.
\bibitem{pik} 
Pinkham, H.C., 
Automorphisms of cusps and Inoue-Hirzebruch surfaces, 
Compositio Math. {\bf 52}  
(1984), 299--313. 
\bibitem{st} 
Steenbrink, J.H.M., 
Mixed Hodge structure on the vanishing cohomology, 
ed. Holm, Nordic Summer School, NAVF, Symposium in Math., Oslo, 1976, 535--563. 
\end{thebibliography}}
\end{document}